\documentclass[11pt,reqno]{amsart}

\usepackage{amsmath, amsfonts, amsthm, amssymb}

\newtheorem{Theorem}{Theorem}[section]
\newtheorem{Lemma}{Lemma}[section]
\newtheorem{Proposition}{Proposition}[section]

\theoremstyle{definition}
\newtheorem{Definition}{Definition}[section]

\theoremstyle{remark}
\newtheorem{Remark}{Remark}[section]

\numberwithin{equation}{section}

\renewcommand{\u}{{\bf u}}

\newcommand{\R}{{\mathbb R}}
\newcommand{\Dv}{{\rm div}}

\newcommand{\dl}{\delta}

\def\f{\frac}
\renewcommand{\O}{\Omega}

\def\D{\Delta }

\def\hf1{^\f{1}{1-\xi^2}}

\def\be{\begin{equation}}
\def\en{\end{equation}}
\def\bs{\begin{split}}
\def\es{\end{split}}

\renewcommand{\d}{{\bf d}}
\newcommand{\F}{{\mathtt F}}
\renewcommand{\v}{{\bf v}}
\renewcommand{\a}{\alpha}

\author{Xianpeng Hu and Dehua Wang}
\address{Department of Mathematics, University of Pittsburgh,
                           Pittsburgh, PA 15260.}
\email{xih15@pitt.edu}
\address{Department of Mathematics, University of Pittsburgh,
                           Pittsburgh, PA 15260.}
\email{dwang@math.pitt.edu}

\title[Global Solution to the  Flow of Liquid Crystals]
{Global Solution to the Three-Dimensional Incompressible Flow of Liquid Crystals}

\keywords{Liquid crystal, strong solution, weak solution, existence, uniqueness}
 \subjclass[2000]{35A05, 76A10, 76D03.}
\date{\today}

\begin{document}

\begin{abstract}
The equations for the three-dimensional incompressible flow of liquid crystals are considered in a smooth bounded domain.
The existence and uniqueness of the global  strong solution with small initial data are established.
It is also proved that  when the strong solution exists, all the global weak solutions
constructed in \cite{LL} must be equal to the unique strong solution.
\end{abstract}
\maketitle

\section{Introduction}

Liquid crystals are substances that exhibit a phase of matter that has properties between those of a conventional liquid, and those of a solid crystal. For instance, a liquid crystal  may flow like a liquid, but its molecules may be oriented in a crystal-like way. There are many different types of liquid crystal phases, which can be distinguished based on their different optical properties.
The various liquid crystal phases can be characterized by the type of ordering that is present. One can distinguish positional order and orientational order, and moreover order can be either short-range or long-range.
Liquid crystals may have an isotropic phase at high temperature, or anisotropic orientational structure at lower temperature.
The diverse phases of liquid crystals have wide applications from the liquid crystal display to biology (In particular, biological membranes and cell membranes are a form of liquid crystal).
In the 1960s, the theoretical physicist P.-G. de Gennes found fascinating analogies between liquid crystals and superconductors as well as magnetic materials, which was rewarded with the Nobel Prize in Physics in 1991.
One of the most common liquid crystal phases is the nematic, where the molecules have no positional order, but they have long-range orientational order. For more details of physics, we refer the readers to the two books of
de Gennes-Prost \cite{Gen} and  Chandrasekhar \cite{Chan}.

The three-dimensional flow of nematic liquid crystals can be governed by the following system of partial differential equations
(\cite{Gen, Lin2,Lin1,LL}):
%
%
%
%
\begin{subequations}\label{e2}
\begin{align}
&\f{\partial\u}{\partial t}+\u\cdot\nabla\u-\mu\D\u+\nabla P
       =-\lambda\Dv\left(\nabla \d\odot\nabla\d\right),\label{e21}\\  
&\f{\partial\d}{\partial t}+\u\cdot\nabla\d=\gamma\left(\D\d-f(\d)\right),\label{e22}\\
&\Dv\u=0.\label{e23}
\end{align}
\end{subequations}
where $\u\in\R^3$ denotes the vector field, $\d\in\R^3$ the director field for the averaged macroscopic molecular orientations,
 $P\in\R$ the pressure arising from the incompressibility;  and they all depend on the spatial variable
 $x=(x_1,x_2,x_3)\in\R^3$ and the time variable $t>0$. 
The positive constants $\mu, \lambda, \gamma$ stand for viscosity, the competition between kinetic energy and potential energy, and microscopic elastic relaxation time or the Debroah number for the molecular orientation field, respectively. We set these three constants to be one since their sizes do not play any role in our analysis.
The symbol $\nabla\d\odot\nabla\d$ denotes a matrix whose $ij$-th entry is $<\partial_{x_i} \d, \partial_{x_j}\d>$, and it is easy to see that $$\nabla\d\odot\nabla\d=(\nabla \d)^\top\nabla\d,$$
where $(\nabla \d)^\top$ denotes the transpose of the $3\times 3$ matric $\nabla \d$.
In \eqref{e2}, $f(\d)$ is the penalty function which will be assumed to be zero as in \cite{LL} for the three-dimensional problem.
The system \eqref{e2} is a simplified version,  but still retains most of the essential features,   of the Ericksen-Leslie equations (\cite{Eri, Eri2, HK, HKL, Leslie1, Leslie2}) for the hydrodynamics of nematic liquid crystals; see \cite{LL,LW,SL} for more discussions on the relations of the two models.
Both the Ericksen-Leslie system and the simplified one \eqref{e2} describe the time evolution of liquid crystal materials under the influence of both the velocity field $\u$ and the director field $\d$. In many situations, the flow velocity field does disturb the alignment of the molecule, and a change in the alignment will induce velocity.

We consider the initial-boundary value problem of system \eqref{e2} in a bounded domain $\O\subset\R^3$
with $C^3$ boundary under the initial-boundary conditions:
\begin{equation}\label{ic}
\d|_{t=0}=\d_0,\quad \u|_{t=0}=\u_0,
\end{equation}
and
\begin{equation}\label{bc}
\u|_{\partial\O}=0,\quad\d|_{\partial\O}=\d_0,
\end{equation}
with $\Dv \u_0=0$ in $\O$, and $\d_0\in C^1(\overline{\O})$ satisfying
$\nabla\d_0=0$ on the boundary $\partial\O$.
We introduce an $3\times 3$ matrix
\begin{equation}\label{1}
\F=\nabla\d,
\end{equation}
and take the gradient of \eqref{e22} to rewrite  \eqref{e2}, with $f(\d)=0$ and $\mu=\lambda=\gamma=1$,  as:
\begin{subequations}\label{e3}
\begin{align}
&\f{\partial\u}{\partial t}+\u\cdot\nabla\u-\D\u+\nabla
P=-\Dv(\F^\top\F),\label{e31}\\
&\f{\partial\F}{\partial
t}+\u\cdot\nabla\F+\F\nabla\u=\D\F,\label{e32}\\
&\Dv\u=0,\label{e33}
\end{align}
\end{subequations}
where  we used, for all $i,j,k=1, 2, 3,$
$$\f{\partial}{\partial x_k}\left(\u_j\f{\partial \d_i}{\partial
x_j}\right)=\f{\partial \u_j}{\partial x_k}\f{\partial
\d_i}{\partial x_j}+\u_j \f{\partial}{\partial
x_j}\left(\f{\partial \d_i}{\partial
x_k}\right)=(\F\nabla\u+\u\cdot\nabla \F)_{ik}.$$
Notice that \eqref{e31} is the incompressible Navier-Stokes
equation with the source term, $-\Dv(\F^\top\F)$, while
\eqref{e32} is a parabolic equation  of $\F$. The initial-boundary
conditions \eqref{ic} and \eqref{bc} become
\begin{equation}\label{ic2}
\u|_{t=0}=\u_0,\quad \F|_{t=0}=\F_0:=\nabla \d_0,
\end{equation}
and
\begin{equation}\label{bc2}
\u|_{\partial\O}=0,\quad\F|_{\partial\O}=0.
\end{equation}
There have been some studies on system \eqref{e2}.  In Lin-Liu
\cite{LL}, the global existence of weak solutions with large
initial data was proved under the condition that the orientational
configuration $\d(x,t)$ belongs to $H^2$, and the global existence
of classical solutions was also obtained if the coefficient $\mu$
is large enough in three dimensional spaces. The similar results
were obtained also in \cite{SL} for a different but similar model.
When weak solutions are discussed, the regularity of the weak
solution was investigated in \cite{LL2} (and also \cite{HKL}).

In this paper, we are interested in strong solutions of \eqref{e3}
in the  Sobolev space $W^{2,q}(\O)$ with $q>3$. It is
worthy of pointing out that if $\F$ belongs to $W^{2,q}(\O)$, it
is equivalent to saying that $\d$ should be in $W^{3,q}(\O)$
according to \eqref{1}. By a $\textit{Strong Solution}$, we means a
triplet $(\u, \F, P)$ satisfying \eqref{e3} almost everywhere with
the initial condition \eqref{ic2} and the boundary condition \eqref{bc2}.
Our strategy to consider \eqref{e3} in $W^{2,q}(\O)$ is to linearize \eqref{e3} as
\begin{subequations}\label{e4}
\begin{align}
&\f{\partial\u}{\partial t}-\D\u+\nabla
P=-\v\cdot\nabla\v-\Dv(G^\top G),\label{e41}\\
&\f{\partial\F}{\partial
t}-\D\F=-\v\cdot\nabla G-G\nabla\v,\label{e42}\\
&\Dv\u=0,  \label{e43}
\end{align}
\end{subequations}
for some given $\v\in\R^3$ and $G\in M^{3\times 3}$. One of the
motivations of making such an linearization is that we can use the
maximal regularity of Stokes equations (\cite{D}) and the
parabolic equation (\cite{HA}). We first use an iteration method
to establish the local existence and uniqueness of strong solution
with general initial data. Then we prove the global existence by
establishing some global estimates under the condition that the
initial data is small in some sense. The global weak solution was
obtained in Lin-Liu \cite{LL}, but the uniqueness is still an open
problem. We shall prove that when the strong solution exists, all
the global weak solutions constructed in \cite{LL} must be equal
to the unique strong solution, which is called the weak-strong
uniqueness. Similar results were obtained by Danchin \cite{D} for
the density-dependent incompressible Navier-Stokes equations. We
shall establish our results in the spirit of \cite{D}, while
developing new estimates for the director field $\d$.

The rest of the  paper is organized as follows. In Section 2, we state
our main results on local and global existence of strong solution, as well as the weak-strong uniqueness.
In Section 3, we recall  the maximal regularity for Stokes equations and the parabolic equation, and also some $L^\infty$ estimates.
In Section 4, we give the proof of the local existence. In Section 5, we prove the global existence.
Finally in Section 6, we  show the weak-strong uniqueness.
\bigskip

\section{Main Results}

In this section, we state our main results.
If $k>0$ is an integer and $p\ge 1$, we denote by $W^{k,p}$ the set of functions in $L^p(\O)$ whose derivatives of up to order
$k$ belong to $L^p(\O)$.
For $T>0$ and a function space $X$, denote by $L^p(0,T; X)$ the set of Bochner measurable X-valued time dependent
functions $f$ such that $t\rightarrow \|f\|_{X}$ belongs to $L^p(0,T)$.
Let us  define the functional spaces in which the existence of solutions
is going to be obtained:

\begin{Definition}\label{df1}
For $T>0$ and $1<p, q<\infty$, we denote by $M^{p,q}_T$ the set of
triplets $(\u,\F, P)$ such that
$$\u\in C([0,T]; D_{A_q}^{1-\f{1}{p},p})\cap L^p(0,T;
W^{2,q}(\O)\cap W^{1,q}_0(\O)), \quad\partial_t\u\in L^p(0,T;
L^q),\quad\Dv\u=0,$$
$$\F\in C([0,T]; B_{q,p}^{2(1-\f{1}{p})}\cap L^p(0,T;W^{2,q}(\O)), \quad \partial_t\F\in L^p(0,T; L^q(\O)),$$
and
$$P\in L^p(0,T; W^{1,q}(\O)), \quad \int_\O P dx=0.$$
The corresponding norm is denoted by $\|\cdot\|_{M^{p,q}_T}$.
\end{Definition}

In the above definition, the space $D_{A_q}^{1-\f{1}{p},p}$ stands for some fractional
domain of the Stokes operator in $L^q$ (cf. Section 2.3 in
\cite{D}). Roughly, the vector-fields of $D_{A_q}^{1-\f{1}{p},p}$
are vectors which have  $2-\f{2}{p}$ derivatives   in $L^q$, are
divergence-free, and vanish on $\partial\O$. The Besov space (for
definition, see \cite{BL}) $B_{q,p}^{2(1-\f{1}{p})}$ can be
regarded as the interpolation space between $L^q$ and $W^{2,q}$,
that is, $$B_{q,p}^{2(1-\f{1}{p})}=(L^q, W^{2,q})_{1-\f{1}{p},p}.$$
We note that, from Proposition 2.5 in \cite{D},
\begin{equation}\label{I1}
D_{A_q}^{1-\f{1}{p},p}\hookrightarrow B^{2(1-\f{1}{p})}_{q,p}\cap
L^q(\O) .
\end{equation}

The local existence will be shown by using an iterative method, and
if the initial data is sufficiently small in some suitable
function spaces, the solution is indeed global in time. More
precisely, our existence results read:

\begin{Theorem}\label{T1}
Let $\O$ be a bounded domain in $\R^3$ with $C^3$ boundary. Assume
$1\le p,q\le\infty$ with $\f{2}{p}(1-\f{3}{q})\in (0,1)$ and
$\u_0\in D_{A_q}^{1-\f{1}{p},p}, \, \F_0\in
B_{q,p}^{2(1-\f{1}{p})}\cap L^q$. Then,
\begin{enumerate}
\item There exists a $T_0>0$,  such that, system \eqref{e3} with the initial-boundary conditions \eqref{ic2}-\eqref{bc2} has a unique local strong solution $(\u, \F, P)\in M^{p,q}_{T_0}$ in $\O\times (0, T_0)$;
\item Moreover,  there exists a $\delta_0>0$, such that, if the initial data  satisfies
$$\|\u_0\|_{D_{A_q}^{1-\f{1}{p},p}}\le\delta_0, \quad
\|\F_0\|_{B_{q,p}^{2(1-\f{1}{p})}\cap L^q}\le\delta_0,$$ then
\eqref{ic2}-\eqref{bc2} has a unique global strong solution $(\u,
\F, P)\in M^{p,q}_{T}$ in $\O\times (0, T)$ for all  $T>0$.
\end{enumerate}
\end{Theorem}

\begin{Remark}   The above Theorem gives us the global strong solution near $\u=0, \, \F=0$.
The similar argument to the proof of Theorem \ref{T1} below will also enable us to show the global existence of strong solution to \eqref{e3} near the equilibrium state: $\u=0$, $\F=I$ (the $3\times 3$ identity matrix).
\end{Remark}

According to Lin-Liu \cite{LL}, for the given initial-boundary conditions
\eqref{ic2} and \eqref{bc2}, there exists at least a $\textit{Weak
Solution}$ to \eqref{e3}. But its uniqueness is still an open question. More
precisely, a triplet $(v, E, \Pi
)$ is called  a weak solution to \eqref{e3} with \eqref{ic2} and \eqref{bc2}
in $\O\times(0,T)$ if $(v, E, \Pi
)$ satisfies the system
\eqref{e3} in the sense of distributions, i.e, for all $\psi\in (C^\infty_0(\O\times(0,T)))^3$ with $\Dv\psi=0$ and
$\phi\in(C^\infty_0(\O\times(0,T)))^9$, we have
\begin{equation*}
\begin{split}
\int_0^T\!\!\!\!\int_\O v\partial_t\psi \,dxdt+\int_0^T\!\!\!\!\int_\O  v \otimes v
:\nabla\psi \,dxdt-\int_0^T\!\!\!\!\int_\O \nabla v:\nabla\psi
\,dxdt=-\int_0^T\!\!\!\!\int_\O  E^\top E:\nabla\psi \,dxdt;
\end{split}
\end{equation*}
and
\begin{equation*}
\begin{split}
\int_0^T E:\partial_t\phi \,dxdt-\int_0^T\!\!\!\!\int_\O  \u\cdot\nabla
E:\phi \,dxdt-\int_0^T\!\!\!\!\int_\O  E\nabla\u:\phi \,dxdt=\int_0^T\!\!\!\!\int_\O
\nabla E:\nabla\phi \,dxdt,
\end{split}
\end{equation*}
with the energy inequality:
\begin{equation*}
\begin{split}
\int_\O(|v(t)|^2+|E(t)|^2)dx+\int_0^t\int_\O(|\nabla v |^2+|\nabla
E |^2)dxds\le \int_\O(|v_0|^2+|E_0|^2)dx.
\end{split}
\end{equation*}
In this weak formulation, the pressure $\Pi
$ can be determined as in the Navier-Stokes equations, see Galdi \cite{Galdi}.
We state here the existence of weak solutions in Theorem A of \cite{LL}:

\begin{Proposition} \label{P1}
Assume that $\u_0\in L^2$ and $\F_0\in L^2$. Then the system
\eqref{e3} with the initial condition \eqref{ic2} and the boundary
condition \eqref{bc2} has a global weak solution $(v, E, \Pi)$ such that
$$v\in L^2(0,T; H^1)\cap L^\infty(0, T; L^2),$$
and
$$E\in L^2(0,T; H^1)\cap L^\infty(0,T; L^2),$$
for all $T\in (0,\infty)$.
\end{Proposition}

For the same initial-boundary conditions, the relation between its
weak solution and its strong solution can be formulated as:

\begin{Theorem}\label{T2}
Assume that $\u_0\in D_{A_q}^{1-\f{1}{p},p}$ and $\F_0\in
B_{q,p}^{2(1-\f{1}{p})}\cap L^q$. Then its corresponding weak
solution to \eqref{e3} with \eqref{ic2} and \eqref{bc2} is unique
and indeed is equal to its unique strong solution.
\end{Theorem}

Usually, we call this kind of uniqueness as $\textit{Weak-Strong
Uniqueness}$. For the similar results on the compressible
Navier-Stokes equations, we refer the  readers to \cite{EB, L}.

\bigskip

\section{Maximal Regularity}

In this section, we recall the maximal regularities for the parabolic operator and the Stokes operator, as well as some $L^\infty$ estimates.
For $T>0$, $1<p,q<\infty$, denote
$$\mathcal{W}(0,T):= W^{1,p}(0,T; (L^q(\O))^3)\cap L^p(0,T;(W^{2,q}(\O))^3).$$
Throughout this paper, $C$ stands for a generic positive constant.

We first recall the maximal regularity for the parabolic operator (cf.
Theorem 4.10.7 and Remark 4.10.9 in \cite{HA}):
\begin{Theorem}\label{T3}
Given $1<p<\infty$, $\omega_0\in B_{q,p}^{2(1-\f{1}{p})}$ and
$f\in L^p(0,T; L^q(\R^3)^3)$, the Cauchy problem
$$\f{d\omega}{dt}-\D\omega=f,\quad t\in(0,T), \quad \omega(0)=\omega_0,$$
has a unique solution $\omega\in \mathcal{W}(0,T)$, and
$$\|\omega\|_{\mathcal{W}(0,T)}\le C\left(\|f\|_{L^p(0,T; L^q(\R^3))}+\|\omega_0\|_{B_{q,p}^{2(1-\f{1}{p})}}\right),$$
where $C$ is independent of $\omega_0$, $f$ and $T$. Moreover,
there exists a positive constant $c_0$ independent of $f$ and $T$
such that
$$\|\omega\|_{\mathcal{W}(0,T)}\ge c_0\sup_{t\in(0,T)}\|\omega(t)\|_{B_{q,p}^{2(1-\f{1}{p})}}.$$
\end{Theorem}

Now we recall the maximal regularity  for the Stokes equations (cf. Theorem 3.2 in \cite{D}):

\begin{Theorem}\label{T4}
Let $\O$ be a bounded domain with $C^3$ boundary in
$\R^3$ and $1<p,q<\infty$. Assume that $\u_0\in
D_{A_q}^{1-\f{1}{p},p}$ and $f\in L^p(\R^{+}; L^q)$. Then the
system
\begin{equation*}
\begin{cases}
\partial_t\u-\D\u+\nabla P=f,\quad \int_\O Pdx=0,\\
\Dv\u=0,\quad \u|_{\partial\O}=0,\\
\u|_{t=0}=\u_0,
\end{cases}
\end{equation*}
has a unique solution $(\u, P)$ satisfying the following
inequality for all $T>0$:
\begin{equation}\label{41}
\begin{split}
\|\u(T)\|_{D_{A_q}^{1-\f{1}{p},p}}+&\left(\int_0^T\left\|\left(\nabla
P, \D\u,
\partial_t\u\right)\right\|^p_{L^q}dt\right)^{\f{1}{p}}\\
&\le
C\left(\|\u_0\|_{D_{A_q}^{1-\f{1}{p},p}}+\left(\int_0^T\|f(t)\|_{L^q}^pdt\right)^{\f{1}{p}}\right)
\end{split}
\end{equation}
with $C=C(q, p, \O)$.
\end{Theorem}

\begin{Remark}
We notice that  \eqref{41} does not include the estimate for $\|\u\|_{L^p(0,T; L^q)}$.
Indeed,  thanks to $\u|_{\partial\O}=0$,  Poincare's inequality,
and the fact $\int_\O\nabla\u dx=0$, we have
$$\|\u\|_{W^{2,q}}\le C\|\D\u\|_{L^q},$$ and then \eqref{41} can
be rewritten as
\begin{equation}\label{41b}
\begin{split}
\|\u(T)\|_{D_{A_q}^{1-\f{1}{p},p}}+&\left(\int_0^T\left\|\left(\nabla
P, \u, \D\u,
\partial_t\u\right)\right\|^p_{L^q}dt\right)^{\f{1}{p}}\\
&\le
C\left(\|\u_0\|_{D_{A_q}^{1-\f{1}{p},p}}+\left(\int_0^T\|f(t)\|_{L^q}^pdt\right)^{\f{1}{p}}\right).
\end{split}
\end{equation}
\end{Remark}

We have the $L^\infty$ estimate in the spatial variable as follows (cf. Lemma 4.1 in \cite{D}).

\begin{Lemma}\label{l1}
Let $1<p,q,r,s<\infty$ satisfy
$$0<\f{p}{2}-\f{3p}{2r}<1, \quad \f{1}{s}=\f{1}{r}+\f{1}{q}.$$
 Then the  following inequalities hold:
$$\|\nabla f\|_{L^p(0,T; L^\infty)}\le
CT^{\f{1}{2}-\f{3}{2r}}\|f\|^{1-\theta}_{L^\infty(0,T;
D_{A_r}^{1-\f{1}{p},p})}\|f\|^\theta_{L^p(0,T; W^{2,r})},$$
$$\|\nabla f\|_{L^p(0,T; L^q)}\le
CT^{\f{1}{2}-\f{3}{2r}}\|f\|^{1-\theta}_{L^\infty(0,T; D_{A_s}^{1-\f{1}{p},p})}\|f\|^\theta_{L^p(0,T; W^{2,s})},$$ for
some constant $C$ depending only on $\O, p, q$ and
$$\f{1-\theta}{p}=\f{1}{2}-\f{3}{2r}.$$
\end{Lemma}

Similarly, we have,
\begin{Lemma}\label{l2}
Let $1<p,q<\infty$ satisfy $0<\f{p}{2}-\f{3p}{2q}<1$. Then one has,
$$\|\nabla f\|_{L^p(0,T; L^\infty)}\le
CT^{\f{1}{2}-\f{3}{2q}}\|f\|^{1-\theta}_{L^\infty(0,T; B_{q,p}^{2(1-\f{1}{p}),p})}
\|f\|^\theta_{L^p(0,T; W^{2,q})},$$
for some constant $C$ depending only on $\O, p, q$ and
$$\f{1-\theta}{p}=\f{1}{2}-\f{3}{2q}.$$
\end{Lemma}
\begin{proof}
First, we notice that
$$(B^{1-\f{2}{p}-\f{3}{q}}_{\infty,\infty}, B^{1-\f{3}{q}}_{\infty,\infty})_{\theta,1}=B^0_{\infty,1}$$
with $$\f{1-\theta}{p}=\f{1}{2}-\f{3}{2q},$$ see Theorem 6.4.5 in
\cite{BL}. Also the imbedding $B^0_{\infty,1}\hookrightarrow
L^\infty$ is true due to Theorem 6.2.4 in \cite{BL}. Hence, one
has
\begin{equation}\label{5}
\|\nabla f\|_{L^\infty}\le C\|\nabla f\|_{B^0_{\infty,1}}\le
C\|\nabla f\|^{\theta}_{B^{1-\f{3}{q}}_{\infty,\infty}}\|\nabla
f\|^{1-\theta}_{B^{1-\f{2}{p}-\f{3}{q}}_{\infty,\infty}}.
\end{equation}
We remark that $$B_{q,p}^{2(1-\f{1}{p})}\hookrightarrow
B_{\infty,\infty}^{2-\f{2}{p}-\f{3}{q}}\hookrightarrow
B_{\infty,\infty}^{1-\f{2}{p}-\f{3}{q}},\quad W^{1,q}\hookrightarrow B_{\infty,\infty}^{1-\f{3}{q}},$$
see Theorem 6.2.4 and Theorem 6.5.1 in \cite{BL}. Hence, according to
\eqref{5}, one deduce that
\begin{equation*}
\begin{split}
\|\nabla f\|_{L^p(0,T; L^\infty)}&\le C\left(\int_0^T\|\nabla
f\|^{p\theta}_{B^{1-\f{3}{q}}_{\infty,\infty}}\|\nabla
f\|^{p(1-\theta)}_{B^{1-\f{2}{p}-\f{3}{q}}_{\infty,\infty}}dt\right)^{\f{1}{p}}\\
&\le
C\left(\int_0^T\|f\|^{p\theta}_{W^{2,q}}\|f\|^{p(1-\theta)}_{B_{q,p}^{2(1-\f{1}{p}),p}}dt\right)^{\f{1}{p}}\\
&\le CT^{\f{1}{2}-\f{3}{2q}}\|f\|^{1-\theta}_{L^\infty(0,T; B_{q,p}^{2(1-\f{1}{p}),p})}\|f\|^\theta_{L^p(0,T; W^{2,q})}.
\end{split}
\end{equation*}
\end{proof}

\begin{Lemma}\label{l3}
For $f\in L^p(0,T, L^q)$ and $\partial_t f\in L^p(0,T; L^q)$ with
$f(0)\in L^q$, we have, for all $t\in[0,T]$,
\begin{equation}\label{41c}
\|f\|_{L^\infty(0,t;L^q)}\le C\left(\|f_0\|_{L^q}+\|f\|_{L^p(0,t; L^q)}+\|\partial_t f\|_{L^p(0,t; L^q)}\right),
\end{equation}
for some positive constant $C$ independent of $T$ and $f$.
\end{Lemma}
\begin{proof}
Indeed, we have
\begin{equation*}
\begin{split}
\|f(t)\|_{L^q}^p&=\|f_0\|_{L^q}^p+\int_0^t\f{d}{ds}\|f(s)\|_{L^q}^p ds\\
&=\|f_0\|^p_{L^q}+\f{p}{q}\int_0^t \left(\|f(t)\|_{L^q}^{p-q}\int_{\R^3} |f(s)|^{q-2}f(s)\partial_sf(t)dx\right)dt\\
&\le\|f_0\|^p_{L^q}+\f{p}{q}\int_0^t\|f(s)\|_{L^q}^{p-1}\|\partial_s f\|_{L^q}ds\\
&\le
\|f_0\|^p_{L^q}+\f{p}{q}\left(\int_0^t\|f\|_{L^q}^pds\right)^{\f{p-1}{p}}\left(\int_0^t\|\partial_s
f\| _{L^q}^pds\right)^{\f{1}{p}},
\end{split}
\end{equation*}
and consequently, \eqref{41c} follows from H\"older's inequality.
\end{proof}

\bigskip

\section{Local Existence}

In this section, we prove the local existence and uniqueness of strong solution in Theorem \ref{T1}. The proof will be divided into several steps, including constructing the approximate solution by iteration, obtaining the uniform estimate, showing the convergence, consistency, and uniqueness.

\subsection{Construction of approximate solutions}
We initialize the construction of approximate solutions by
setting $F^0:= \F_0$ and $\u^0:=\u_0$. For given
$(\u^n, \F^n)$,  the Stokes equations \eqref{e41} and the parabolic
equation \eqref{e42} enable us to define $(\u^{n+1}, \F^{n+1}, P^{n+1})$
as the global solution of
\begin{subequations}\label{e5}
\begin{align}
&\f{\partial\u^{n+1}}{\partial t}-\D\u^{n+1}+\nabla
P^{n+1}=-\u^{n}\cdot\nabla\u^{n}-\Dv({\F^{n}}^\top\F^{n}),\label{e51}\\
&\f{\partial\F^{n+1}}{\partial
t}-\D\F^{n+1}=-\u^{n}\cdot\nabla \F^{n}-\F^{n}\nabla\u^{n},\label{e52}\\
&\Dv\u^{n+1}=0,\label{e53}
\end{align}
\end{subequations}
with the initial-boundary conditions:
$$\u^{n+1}|_{t=0}=\u_0,\quad\F^{n+1}|_{t=0}=\F_0,$$
$$\u^{n+1}|_{\partial\O}=0,\quad \F^{n+1}|_{\partial\O}=0,$$
and
$$\int_{\O}P^{n+1}dx=0.$$

According to Theorem \ref{T3} and Theorem \ref{T4}, an argument by
induction yields a sequence $\{(\u^n, \F^n, P^n)\}_{n\in {\mathbb N}} \subset M^{p,q}_T$ for all positive $T$.

\subsection{ Uniform estimate for some small fixed time $T$}   We aim
at finding a positive time $T$ independent of $n$ for which
$\{(\u^n, \F^n, P^n)\}_{n\in {\mathbb N}}
$ is uniformly bounded in the space
$M^{p,q}_T$. Indeed, applying Theorem \ref{T3} and Theorem
\ref{T4}, we obtain
\begin{equation}\label{3}
\begin{split}
&\|\u^{n+1}(T)\|_{D_{A_q}^{1-\f{1}{p},p}}+\left(\int_0^T\left\|\left(\nabla
P^{n+1}, \u^{n+1}, \D\u^{n+1},
\partial_t\u^{n+1}\right)\right\|^p_{L^q}dt\right)^{\f{1}{p}}\\
&\le
C\left(\|\u_0\|_{D_{A_q}^{1-\f{1}{p},p}}+\left(\int_0^T\|\u^n\cdot\nabla\u^n+\Dv({\F^n}^\top\F^{n})\|_{L^q}^pdt\right)^{\f{1}{p}}\right),
\end{split}
\end{equation}
and
\begin{equation}\label{4}
\begin{split}
&\|\F^{n+1}(T)\|_{B_{q,p}^{2(1-\f{1}{p})}}+\|\F^{n+1}\|_{\mathcal{W}(0,T)}\\
&\le
C\left(\|F_0\|_{B_{q,p}^{2(1-\f{1}{p})}}+\|\F^n\nabla\u^n+\u^n\cdot\nabla\F^n\|_{L^p(0,T;
L^q(\O))}\right).
\end{split}
\end{equation}

Now define
\begin{equation*}
\begin{split}
U^n(t):=&\|\u^n(t)\|_{L^\infty(0,t; D_{A_q}^{1-\f{1}{p},p})}+\|\u^n\|_{L^p(0,t;
W^{2,q})}+\|\partial_t\u^n\|_{L^p(0,t; L^q)}\\
&+\|\F^n(t)\|_{L^\infty(0,t;
B_{q,p}^{2(1-\f{1}{p})})}+\|\F^n\|_{\mathcal{W}(0,t)},
\end{split}
\end{equation*}
 and
$$U^0\:=\|\u_0\|_{D_{A_q}^{1-\f{1}{p},p}}+\|F_0\|_{B_{q,p}^{2(1-\f{1}{p})}\cap L^q}.$$

Hence, from \eqref{3} and \eqref{4}, one has, using Lemmas \ref{l1}-\ref{l3},
\begin{equation}\label{6}
\begin{split}
U^{n+1}(t)&\le C\Big(U^0+\|\F^n\|_{L^\infty(0,t;
L^q)}\|\nabla\F^n\|_{L^p(0,t; L^\infty)}+\|\u^n\|_{L^\infty(0,t;
L^q)}\|\nabla\u^n\|_{L^p(0,t; L^\infty)}\\
&\qquad\quad+\|\u^n\|_{L^\infty(0,t;
L^q)}\|\nabla\F^n\|_{L^p(0,t; L^\infty)}+\|\F^n\|_{L^\infty(0,t;
L^q)}\|\nabla\u^n\|_{L^p(0,T; L^\infty)}\Big)\\
&\le C\left(U^0+t^{\f{1}{2}-\f{3}{2q}}(U^0+U^n(t))U^n(t)\right).
\end{split}
\end{equation}
Hence, if we assume that $U^n(t)\le 4CU^0$ on $[0,T_0]$ with
\begin{equation}\label{TTT1}
0<T_0\le \left(\f{3}{4C(4C+1)U^0}\right)^{\f{2q}{q-3}},
\end{equation}
 then a direct computation yields
\begin{equation*}
U^{n+1}(t)\le 4CU^0,\quad\textrm{on}\quad [0,T_0].
\end{equation*}

Coming back to \eqref{3}, \eqref{4}, and \eqref{6}, we conclude
that the sequence $\{(\u^n, \F^n, P^n)\}_{n=1}^\infty$ is
uniformly bounded in $M^{p,q}_{T_0}$. More precisely, we have
\begin{Lemma}
For all $t\in [0,T_0]$ with $T_0$ satisfying \eqref{TTT1},
\begin{equation}\label{7}
U^n(t)\le 4CU^0.
\end{equation}
\end{Lemma}

\subsection{Convergence of the approximate sequence}
We now prove
\begin{Lemma}
 $\{(\u^n, \F^n, P^n)\}_{n=1}^\infty$ is a
Cauchy sequence and thus converges in  $M_{T_0}^{p,q}$.
\end{Lemma}
\begin{proof}
Let
$$\dl\u^n:=\u^{n+1}-\u^n,\quad \dl P^n:=P^{n+1}-P^n,\quad \dl\F^n:=\F^{n+1}-\F^n.$$
 Define
\begin{equation}\label{14}
\begin{split}
\dl U^n(t):=&\|\dl\u^n(t)\|_{L^\infty(0,t; D_{A_q}^{1-\f{1}{p},p})}+\|\dl\u^n\|_{L^p(0,t;
W^{2,q})}+\|\partial_t\dl\u^n\|_{L^p(0,t; L^q)}\\
&+\|\dl\F^n(t)\|_{L^\infty(0,t; B_{q,p}^{2(1-\f{1}{p})})}+\|\dl\F^n\|_{\mathcal{W}(0,t)}.
\end{split}
\end{equation}
The triplet $(\dl\u^n, \dl\F^n, \dl P^n)$ satisfies
\begin{equation}\label{e6}
\begin{cases}
\f{\partial\dl\u^{n}}{\partial t}-\D\dl\u^{n}+\nabla \dl P^{n}\\
\qquad\qquad =-\u^{n}\cdot\nabla\u^{n}+\u^{n-1}\cdot\nabla\u^{n-1}-\Dv({\F^{n}}^\top\F^{n})+\Dv({\F^{n-1}}^\top\F^{n-1}),\\
\f{\partial\dl\F^{n}}{\partial
t}-\D\dl\F^{n}=-\u^{n}\cdot\nabla \F^{n}-\F^{n}\nabla\u^{n}+\u^{n-1}\cdot\nabla \F^{n-1}+\F^{n-1}\nabla\u^{n-1},\\
\Dv\u^{n}=0,
\end{cases}
\end{equation}
with
$$\dl\u^n|_{t=0}=\dl\u^n|_{\partial\O}=0,$$
$$\dl\F^n|_{t=0}=\dl\F^n|_{\partial\O}=0,$$
and
$$\int_\O\dl P^ndx=0.$$

Notice that, using Lemma \ref{l1} and Lemma \ref{l2},
\begin{equation}\label{8}
\begin{split}
&\|-\u^{n}\cdot\nabla\u^{n}+\u^{n-1}\cdot\nabla\u^{n-1}\|_{L^p(0,T;
L^q(\O))}\\&=\|\dl\u^{n-1}\cdot\nabla
\u^{n}-\u^{n-1}\cdot\nabla\dl\u^{n-1}\|_{L^p(0,T; L^q)}\\
&\le\|\u^{n-1}\|_{L^\infty(0,T;
L^q)}\|\nabla\dl\u^{n-1}\|_{L^p(0,T;
L^\infty)}+\|\dl\u^{n-1}\|_{L^\infty(0,T;
L^q)}\|\nabla\u^{n}\|_{L^p(0,T; L^\infty)}\\
&\le 4CU^0\left(\|\nabla\dl\u^{n-1}\|_{L^p(0,T;
L^\infty)}+T^{\f{1}{2}-\f{3}{2q}}\|\dl\u^{n-1}\|_{L^\infty(0,T; L^q)}\right),
\end{split}
\end{equation}
\begin{equation}\label{9}
\begin{split}
&\|-\Dv({\F^{n}}^\top\F^{n})+\Dv({\F^{n-1}}^\top\F^{n-1})\|_{L^p(0,T;
L^q(\O))}\\&=\|-\Dv({\dl\F^{n-1}}^\top\F^{n})-\Dv({\F^{n-1}}^\top\dl\F^{n-1})\|_{L^p(0,T; L^q)}\\
&\le\|\F^n\|_{L^\infty(0,T; L^q)}\|\nabla\dl\F^{n-1}\|_{L^p(0,T;
L^\infty)}+\|\dl\F^{n-1}\|_{L^\infty(0,T;
L^q)}\|\nabla\F^{n-1}\|_{L^p(0,T; L^\infty)}\\
&\quad+\|\F^{n-1}\|_{L^\infty(0,T;
L^q)}\|\nabla\dl\F^{n-1}\|_{L^p(0,T;
L^\infty)}+\|\dl\F^{n-1}\|_{L^\infty(0,T;
L^q)}\|\nabla\F^{n}\|_{L^p(0,T; L^\infty)}\\
&\le 4CU^0\left(\|\nabla\dl\F^{n-1}\|_{L^p(0,T;
L^\infty)}+T^{\f{1}{2}-\f{3}{2q}}\|\dl\F^{n-1}\|_{L^\infty(0,T; L^q)}\right),
\end{split}
\end{equation}
\begin{equation}\label{10}
\begin{split}
&\|-\u^n\cdot\nabla\F^n+\u^{n-1}\cdot\nabla\F^{n-1}\|_{L^p(0,T;
L^q)}\\&=\|\u^n\cdot\nabla\dl\F^{n-1}+\dl\u^{n-1}\cdot\nabla\F^{n-1}\|_{L^p(0,T;
L^q)}\\&\le\|\u^n\|_{L^\infty(0,T;
L^q)}\|\nabla\dl\F^{n-1}\|_{L^p(0,T;
L^\infty)}+\|\dl\u^{n-1}\|_{L^\infty(0,T;
L^q)}\|\nabla\F^{n-1}\|_{L^p(0,T; L^\infty)}\\
&\le 4CU^0\left(\|\nabla\dl\F^{n-1}\|_{L^p(0,T;
L^\infty)}+T^{\f{1}{2}-\f{3}{2q}}\|\dl\u ^{n-1}\|_{L^\infty(0,T; L^q)}\right),
\end{split}
\end{equation}
and
\begin{equation}\label{11}
\begin{split}
&\|-\F^n\nabla\u^n+\F^{n-1}\nabla\u^{n-1}\|_{L^p(0,T;
L^q)}\\&=\|\F^n\nabla\dl\u^{n-1}+\dl\F^{n-1}\nabla\u^{n-1}\|_{L^p(0,T;
L^q)}\\&\le\|\F^n\|_{L^\infty(0,T;
L^q)}\|\nabla\dl\u^{n-1}\|_{L^p(0,T;
L^\infty)}+\|\dl\F^{n-1}\|_{L^\infty(0,T;
L^q)}\|\nabla\u^{n-1}\|_{L^p(0,T; L^\infty)}\\
&\le 4CU^0\left(\|\nabla\dl\u^{n-1}\|_{L^p(0,T;
L^\infty)}+T^{\f{1}{2}-\f{3}{2q}}\|\dl\F^{n-1}\|_{L^\infty(0,T; L^q)}\right).
\end{split}
\end{equation}

Applying Theorems \ref{T3}-\ref{T4} with the help of
\eqref{8}-\eqref{11}, one deduce that
\begin{equation}\label{12}
\begin{split}
\dl U^n(t)&\le 8CU^0\Big(\|\nabla\dl\u^{n-1}\|_{L^p(0,t;
L^\infty)}+\|\nabla\dl\F^{n-1}\|_{L^p(0,t;
L^\infty)}\\&\qquad\qquad+t^{\f{1}{2}-\f{3}{2q}}(\|\dl\F^{n-1}\|_{L^\infty(0,t;
L^q)}+\|\dl\u^{n-1}\|_{L^\infty(0,t; L^q)})\Big).
\end{split}
\end{equation}
On the other hand, \eqref{14} implies that, by Lemma \ref{l3},
$$\|\dl\F^{n-1}\|_{L^\infty(0,t;
L^q)}+\|\dl\u^{n-1}\|_{L^\infty(0,t; L^q)}\le \dl U^{n-1}(t),$$
which, combining with \eqref{12}, Lemma \ref{l1} and Lemma \ref{l2}
together, gives
\begin{equation}\label{13}
\dl U^n(t)\le 16CU^0t^{\f{1}{2}-\f{3}{2q}}\dl U^{n-1}(t).
\end{equation}
Thus, if we choose $T_0$ satisfying \eqref{TTT1},  such that,  the condition
$$16CU^0{T_0}^{\f{1}{2}-\f{3}{2q}}\le\f{1}{2}$$ is fulfilled, it is
clear that $\{(\u^n, \F^n, P^n)\}_{n=1}^\infty$ is a Cauchy
sequence in $M_{T_0}^{p,q}$.
\end{proof}

\subsection{The Limit is a solution}
Since $\{(\u^n, \F^n, P^n)\}_{n=1}^\infty$ is a Cauchy sequence in $M_{T_0}^{p,q}$, then it converges.
 Let $(\u,\F, P)\in M_{{T_0}}^{p,q}$ be the limit of the sequence $\{(\u^n, \F^n, P^n)\}_{n=1}^\infty$ in $M_{T_0}^{p,q}$.
 We claim all those nonlinear terms in
\eqref{e5} converge to their corresponding terms in \eqref{e3} in
$L^p(0,{T_0}; L^q)$ . Indeed, using Lemmas \ref{l1} and  \ref{l3}, we have,
\begin{equation*}
\begin{split}
&\|\u^n\cdot\nabla\u^n-\u\cdot\nabla\u\|_{L^p(0,{T_0};
L^q)}\\&=\|(\u^n-\u)\cdot\nabla\u^n+\u\cdot\nabla(\u^n-\u)\|_{L^P(0,{T_0};
L^q)}\\
&\le\|\u^n-\u\|_{L^\infty(0,{T_0}; L^q)}\|\nabla\u^n\|_{L^p(0,{T_0};
L^\infty)}+\|\u\|_{L^\infty(0,{T_0};
L^q)}\|\nabla\u^n-\nabla\u\|_{L^p(0,{T_0}; L^\infty)}\\
&\le C\|\u^n-\u\|_{M_{{T_0}}^{p,q}}{T_0}^{\f{1}{2}-\f{3}{2q}}CU^0
+C\|\u\|_{L^\infty(0,{T_0};
L^q)}{T_0}^{\f{1}{2}-\f{3}{2q}}\|\u^n-\u\|_{M_{{T_0}}^{p,q}}\\&\rightarrow
0,
\end{split}
\end{equation*}
as $n\rightarrow\infty$ due to the convergence of $\u^n$ to $\u$
in $M_{{T_0}}^{p,q}$ and Lemma \ref{l3}. Hence,
$$\u^n\cdot\nabla\u^n\rightarrow\u\cdot\nabla\u,\quad\textrm{in}\quad
L^p(0,{T_0}; L^q).$$
Similarly, we have
$$\Dv(\F^n{\F^n}^\top)\rightarrow\Dv(\F\F^\top),\quad\textrm{in}\quad
L^p(0,{T_0}; L^q);$$
$$\u^n\cdot\nabla\F^n\rightarrow\u\cdot\nabla\F,\quad\textrm{in}\quad
L^p(0,{T_0}; L^q);$$
$$\F^n\nabla\u^n\rightarrow\F\nabla\u,\quad\textrm{in}\quad
L^p(0,{T_0}; L^q).$$
Thus, taking the limit as $n\rightarrow\infty$ in \eqref{e5}, we
conclude that \eqref{e3} holds in $L^p(0,{T_0}; L^q)$, and hence
almost everywhere on $\O\times [0,{T_0}]$.

\subsection{Uniqueness}   Let $(\u_1, \F_1, P_1)$ and $(\u_2,
\F_2, P_2)$ be two solutions to \eqref{e3} with the
initial-boundary conditions \eqref{ic2} and \eqref{bc2}.
Denote $$\dl\u=\u_1-\u_2,\quad \dl\F=\F_1-\F_2,\quad \dl P=P_1-P_2.$$
Note that the triplet $(\dl\u, \dl\F, \dl P)$ satisfies the
following system:
\begin{equation}\label{15}
\begin{cases}
\partial_t\dl\u-\mu\D\dl\u+\nabla\dl
P=-\u_2\cdot\nabla\dl\u-\dl\u\cdot\nabla\u_1+\Dv((\dl\F)^\top\F_1+\F_2^\top\dl\F),\\
\partial_t\dl\F-\D\dl\F=-\u_1\cdot\nabla\dl\F-\dl\u\cdot\nabla\F_2-\F_1\nabla\dl\u-\dl\F\nabla\u_2,\\
\Dv\dl\u=0,
\end{cases}
\end{equation}
with the initial-boundary conditions
$$\dl\u|_{t=0}=\dl\u|_{\partial\O}=0,$$
$$\dl\F|_{t=0}=\dl\F|_{\partial\O}=0,$$
and
$$\int_\O \dl Pdx=0.$$
Define
\begin{equation*}
\begin{split}
X(t):=&\|\dl\u(t)\|_{L^\infty(0,t; D_{A_q}^{1-\f{1}{p},p})}+\|\dl\u\|_{L^p(0,t;
W^{2,q})}+\|\partial_t\dl\u\|_{L^p(0,t; L^q)}\\
&+\|\dl\F(t)\|_{L^\infty(0,t; B_{q,p}^{2(1-\f{1}{p})})}+\|\dl\F\|_{\mathcal{W}(0,t)}.
\end{split}
\end{equation*}
Thus, applying Lemmas \ref{l1} and  \ref{l2} to \eqref{15},
one has, repeating the argument in \eqref{8}-\eqref{11},
\begin{equation*}
\begin{split}
X(t)&\le 4CU^0\Big(\|\nabla\dl\u\|_{L^p(0,t;
L^\infty)}+\|\nabla\dl\F\|_{L^p(0,t;
L^\infty)}\\&\quad+t^{\f{1}{2}-\f{3}{2q}}(\|\dl\F\|_{L^\infty(0,t;
L^q)}+\|\dl\u\|_{L^\infty(0,t; L^q)})\Big)\\
&\le 16CU^0t^{\f{1}{2}-\f{3}{2q}}X(t)\le\f{1}{2}X(t).
\end{split}
\end{equation*}
Hence, $X(t)=0$ for all $t\in [0,{T_0}]$, which guarantee the
uniqueness on the the interval $[0,{T_0}]$.

\bigskip
\section{Global Existence}

In this section, we prove that, if  the initial data is
sufficiently small, the local solution established in the previous
section is indeed global in time. To this end, we first denote
by $T^*$ the maximal time of existence for $(\u, \F, P)$. Define
the function $H(t)$ as
\begin{equation*}
\begin{split}
H(t):=&\|\u\|_{L^p(0,t; W^{2,q})}+\|\partial_t\u\|_{L^p(0,t;
L^q)}+\|\u\|_{L^\infty(0,t;
D_{A_q}^{1-\f{1}{p},p})}\\&\quad+\|P\|_{L^p(0,t;
W^{2,q})}+\|\F\|_{L^\infty(0,t;
B_{q,p}^{2(1-\f{1}{p})})}+\|\F\|_{\mathcal{W}(0,t)},
\end{split}
\end{equation*}
and
$$H_0:=\|\u_0\|_{D_{A_q}^{1-\f{1}{p},p}}+\|\F_0\|_{B_{q,p}^{2(1-\f{1}{p})}\cap L^q}.$$

To extend the local solution, we need to control the maximal time
$T^*$ only in term of the initial data. For this purpose, it is
obvious to observe that $H(t)$ is an increasing and continuous
function in $[0, T^*)$, and for all $t\in [0,T^*)$, we have, using Lemmas \ref{l1} and  \ref{l2},
\begin{equation}\label{18}
\begin{split}
H(t)\le C\Big(H_0+\|\u\cdot\nabla\u\|_{L^p(0,t;
L^q)}+\|\Dv(\F^\top\F)\|_{L^p(0,t; L^q)}
+\|\u\cdot\nabla\F+\F\nabla\u\|_{L^p(0,t; L^q)}\Big).
\end{split}
\end{equation}
On the other hand, Lemmas \ref{l1}-\ref{l3} imply that
\begin{equation}\label{20}
\begin{split}
\|\u\cdot\nabla\u\|_{L^p(0,t; L^q)}&\le \|\u\|_{L^\infty(0,t;
L^q(\O))}\|\nabla\u\|_{L^p(0,t; L^\infty)}\\&\le
C\left(\|u_0\|_{L^q}+H(t)\right)H(t)t^{\f{1}{2}-\f{3}{2q}}\\
&\le C(H_0+H(t))H(t)t^{\f{1}{2}-\f{3}{2q}},
\end{split}
\end{equation}
\begin{equation}\label{21}
\begin{split}
\|\Dv(\F^\top\F)\|_{L^p(0,t; L^q)}&\le C\|F\|_{L^\infty(0,t;
L^q)}\|\nabla\F\|_{L^p(0,t; L^\infty)}\\&\le
C(\|\F_0\|_{L^q}+H(t))H(t)t^{\f{1}{2}-\f{3}{2q}}\\
&\le C(H_0+H(t))H(t)t^{\f{1}{2}-\f{3}{2q}},
\end{split}
\end{equation}
and, similarly, by Lemma \ref{l3},
\begin{equation}\label{22}
\begin{split}
&\|\u\cdot\nabla\F+\F\nabla\u\|_{L^p(0,t; L^q)}\\
&\le \|\u\|_{L^\infty(0,t; L^q)}\|\nabla\F\|_{L^p(0,t;
L^\infty)}+\|\F\|_{L^\infty(0,t; L^q)}\|\nabla\u\|_{L^p(0,t;
L^\infty)}\\
&\le
C(\|\u_0\|_{L^q}+H(t))H(t)t^{\f{1}{2}-\f{3}{2q}}+C(\|\F_0\|_{L^q}+H(t))H(t)t^{\f{1}{2}-\f{3}{2q}}\\
&\le C(H_0+H(t))H(t)t^{\f{1}{2}-\f{3}{2q}}.
\end{split}
\end{equation}
Substituting \eqref{20}-\eqref{22} into \eqref{18}, we get
\begin{equation}\label{23}
H(t)\le C\left(H_0+(H_0+H(t))H(t)t^{\f{1}{2}-\f{3}{2q}}\right).
\end{equation}

Assume that $T$ is the smallest number such that
$$H(T)=4CH_0.$$ This is possible because $H(t)$ is an increasing and
continuous function in time. Then,
$$H(t)<H(T)=4CH_0,\quad\textrm{for all}\quad t\in[0,T),$$ and
from \eqref{23}, we deduce that
$$3\le (H_0+4CH_0)4CT^{\f{1}{2}-\f{3}{2q}}.$$
Hence, we have
$$T^*>T\ge \left(\f{3}{8C(H_0+4CH_0)}\right)^{\f{2q}{q-3}}.$$
This implies that the maximal time of existence will go to
infinity when the initial data  approaches zero.
More precisely,  we can show that,  if the initial
data is sufficiently small, the solution exists globally in time.
To this end, we need some other estimates for the terms on
the right side of \eqref{18}. Indeed, by the imbedding
$$W^{1,q}\hookrightarrow L^\infty,$$
as $q>3$, we have
\begin{equation*}
\begin{split}
\|\u\cdot\nabla\u\|_{L^p(0,t; L^q)}&\le \|\u\|_{L^\infty(0,t;
L^q(\O))}\|\nabla\u\|_{L^p(0,t; L^\infty)}\\&\le
C(\|\u_0\|_{L^q}+H(t))\|\u\|_{L^p(0,t; W^{2,q})}\\
&\le C(H_0+H(t))H(t).
\end{split}
\end{equation*}
Similarly, we have
$$\|\Dv(\F^\top\F)\|_{L^p(0,t; L^q)}\le C(H_0+H(t))H(t),$$
and
$$\|\u\cdot\nabla\F+\F\nabla\u\|_{L^p(0,t; L^q)}\le C(H_0+H(t))H(t).$$
Thus, \eqref{18} turns out to be
\begin{equation}\label{40}
H(t)\le C(H_0+(H_0+H(t))H(t)).
\end{equation}
By the Cauchy-Schwarz inequality, \eqref{40} becomes
\begin{equation}\label{37}
\begin{split}
H(t)&\le C(H_0+H_0^2+2H^2(t)),
\end{split}
\end{equation}
for all $t\in [0,T^*)$. Now we take $H_0$ sufficiently small such that
\begin{equation}\label{38}
H_0+H_0^2\le\dl:=\f{1}{8C^2}.
\end{equation}
Then, under the assumption \eqref{38}, we compute directly from
\eqref{37} and the continuity of $H(t)$ that
\begin{equation}\label{39}
H(t)\le \f{1-\sqrt{1-8C^2(H_0+H_0^2)}}{4C}\le \f{1}{4C},
\end{equation}
for all $t\in [0,T^*)$. In particular, this implies that
$$\|(\u, \F, P)\|_{M^{p,q}_{T^*}}\le\f{1}{4C}<\infty.$$
Hence, according to the local existence in the previous section, we
can extend the solution on $[0,T^*)$ to some larger interval $[0,
T^*+{T_0})$ with ${T_0}>0$. This is impossible since
 $T^*$ is already the maximal time of existence. Hence, when the initial data satisfies \eqref{38},
 the strong solution is indeed global in time.

 The proof of Theorem \ref{T1} is complete.

\bigskip
\section{Weak-Strong Uniqueness}

The purpose of this section is to show $\textit{Weak-Strong Uniqueness}$ in Theorem \ref{T2}.
To this end, we need to obtain first an energy estimate
for the strong solution to the system \eqref{e3}. More precisely,
we have

\begin{Lemma}\label{l4}
Let $p,q$ satisfy the same conditions as Theorem \ref{T1} and
$(\u, \F, P)\in M_{{T_0}}^{p,q}$ be the unique solution to
\eqref{e3} on $\O\times [0, {T_0}]$. Then,  one has,
\begin{equation*}
\int_\O\left(\|\u(t)\|^2+\|\F(t)\|^2\right)dx+\int_0^t\int_\O
\left(\|\nabla\u\|^2+\|\nabla\F\|^2\right)dxds=\int_\O\left(\|\u_0\|^2+\|\F_0\|^2\right)dx.
\end{equation*}
\end{Lemma}

\begin{proof}
Note that $\u\in C([0,{T_0}]; D_{A_q}^{1-\f{1}{p},p})\cap
L^p(0,{T_0}; W^{2,q})$ with $q>3$. Then  $$\u\in C([0,{T_0}];
L^2)\cap L^2(0,{T_0}; H^{1+\a})$$ for some  $\a\geq 0$, since
$$D_{A_q}^{1-\f{1}{p},p}\hookrightarrow B^{2(1-\f{1}{p})}_{q,p}\cap
L^q(\O)\hookrightarrow L^2(\O),$$ Sobolev's embedding
$W^{2,q}(\O)\hookrightarrow H^{2}(\O)$ as $q>3$ and the rest
follows from directly the standard interpolation inequality.
Similarly,
$$\F\in C([0,{T_0}]; L^2)\cap L^2(0,{T_0}; H^{1+\a}).$$

Taking the $L^2$ scalar product in \eqref{e31} with $\u$ and
performing integration by parts, we obtain
\begin{equation}\label{16}
\f{d}{dt}\int_\O|\u|^2dx+\int_\O|\nabla\u|^2dx=\int_\O\F^\top\F:\nabla\u
dx,
\end{equation}
where the notation $A:B$ means the inner product between two
matrix, i.e. $A:B=\sum_{i,j}A_{ij}B_{ij}.$
Similarly, taking the $L^2$ inner product in \eqref{e32} with $\F$
and performing integration by parts, we obtain
\begin{equation}\label{17}
\f{d}{dt}\int_\O|\F|^2dx+\int_\O|\nabla\F|^2dx=-\int_\O
\F\nabla\u:\F dx-\int_\O \F:(\u\cdot\nabla\F)dx,
\end{equation}
where $|\F|^2=\F:\F$ and
$$|\nabla\F|^2=\sum_{i,j,k}\left|\f{\partial\F_{ij}}{\partial x_k}\right|^2.$$

Notice that
$$\int_\O
\F:(\u\cdot\nabla\F)dx=\f{1}{2}\int_\O\u\cdot\nabla|\F|^2dx=-\f{1}{2}\int_\O\Dv\u
|\F|^2dx=0,$$ and, due to $AB:C=A:CB^\top=B:A^\top C$,
$$\int_\O
\F\nabla\u:\F dx=\int_\O\nabla\u:\F^\top\F dx.$$ Hence, adding
\eqref{16} and \eqref{17} together, we have
$$\f{d}{dt}\int_\O(|\u|^2+|\F|^2)dx+\int_\O(|\nabla\u|^2+|\nabla\F|^2)dx=0.$$
Integrating the above equality over time interval $[0,t]$, we
obtain the energy equality of this lemma.
\end{proof}

Now, we recall that for the weak solution $(v,E,\Pi
)$ obtained in \cite{LL},  we have for (almost) all $t\in(0,T)$,
\begin{equation}\label{26}
\begin{split}
\f{1}{2}\int_\O(|v(t)|^2+|E(t)|^2)dx+\int_0^t\int_\O(|\nabla v
|^2+|\nabla E|^2)dxds\le \f{1}{2}\int_\O(|\u_0|^2+|\F_0|^2)dx.
\end{split}
\end{equation}
We remark that, in view of the regularity of $\u$, we deduce
from the weak formulation of \eqref{e3} the following equalities:
\begin{equation}\label{27}
\begin{split}
&\int_\O v\cdot \u dxds+\int_0^t\int_\O\nabla\u:\nabla v dxds\\
&=\int_\O|\u_0|^2+\int_0^t\int_\O E^\top E:\nabla \u dxds
+\int_0^t\int_\O v\cdot\left(\f{\partial \u}{\partial t}+v\cdot\nabla \u\right)dxds,
\end{split}
\end{equation}
and
\begin{equation}\label{28}
\begin{split}
&\int_\O \F: Edx+\int_0^t\int_\O\nabla\F:\nabla Edxds \\
&=\int_\O|\F_0|^2dx-\int_0^t\int_\O v \cdot\nabla E :\F
dxds-\int_0^t\int_\O E\nabla v: \F dxds\\&\quad+\int_0^t\int_\O
E:\f{\partial \F}{\partial t} dxds,
\end{split}
\end{equation}
for a.e. $t\in (0,T)$. Here, we used the identity
$$\int_\O v\cdot\nabla \u\cdot w dx=-\int_\O v\cdot\nabla w\cdot
\u dx,$$ if $\Dv v=0$.

Since $E$ satisfies the equation \eqref{e32}, we substitute \eqref{e32}
into \eqref{28}, and use the following two facts:
$$\int_0^t\int_\O(v\cdot\nabla E:\F+v\cdot\nabla
\F:E)dxds=\int_0^t\int_\O v\cdot\nabla (E:\F) dxds=0,$$ and
$$E\nabla \u:\F+E: \F\nabla \u=\nabla \u:(E^\top \F+\F^\top E),$$
to obtain
\begin{equation}\label{33}
\begin{split}
&\int_\O \F: Edx+2\int_0^t\int_\O\nabla\F:\nabla Edxds \\
&=\int_\O|\F_0|^2dx-\int_0^t\int_\O\nabla \u:(E^\top
\F+\F^\top E)dxds\\&\quad+\int_0^t\int_\O(v-\u)\cdot\nabla \F:E
dxds-\int_0^t\int_\O \F:E\nabla(v-\u) dxds.
\end{split}
\end{equation}

On the other hand, we can write the equation for $\u$ as 
\begin{equation}\label{29}
\begin{split}
\f{\partial \u}{\partial t}+v\cdot\nabla\u-\D \u+\nabla
P=(v-\u)\cdot\nabla \u-\Dv(\F^\top \F).
\end{split}
\end{equation}
Multiplying \eqref{29} by $v$ and integrating over
$\O\times(0,t)$, we get
\begin{equation}\label{30}
\begin{split}
&\int_0^t\int_\O v\cdot\left(\f{\partial \u}{\partial
t}+v\cdot\nabla \u\right)dxds \\
&=-\int_0^t\int_\O\nabla\u:\nabla v
dxds+\int_0^t\int_\O(v-\u)\cdot\nabla \u\cdot v
dxds\\&\quad+\int_0^t\int_\O \F^\top \F:\nabla v dxds.
\end{split}
\end{equation}

Substituting \eqref{30} into \eqref{27}, we obtain
\begin{equation}\label{31}
\begin{split}
&\int_\O\u\cdot v dxds+2\int_0^t\int_\O\nabla\u:\nabla v dxds \\
&=\int_\O|\u_0|^2+\int_0^t\int_\O E^\top E:\nabla \u
dxds\\&\quad+\int_0^t\int_\O(v-\u)\cdot\nabla \u\cdot v
dxds+\int_0^t\int_\O \F^\top \F:\nabla v dxds.,
\end{split}
\end{equation}
Also, according to Lemma \ref{l4}, we have
\begin{equation}\label{32}
\begin{split}
\f{1}{2}\int_\O(|\u|^2+|\F|^2)dx+\int_0^t\int_\O(|\nabla \F
|^2+|\nabla \u|^2)dxds=\f{1}{2}\int_\O(|\u_0|^2+|\F_0|^2)dx.
\end{split}
\end{equation}

Summing \eqref{26}, \eqref{32} and subtracting the sum of
\eqref{33} and \eqref{31}, we obtain for almost all $t\in(0,T)$,
\begin{equation}\label{34}
\begin{split}
&\f{1}{2}\int_\O(|\u(t)-v(t)|^2+|\F(t)-E(t)|^2)dx+\int_0^t\int_\O(|\nabla\u-\nabla
v|^2+|\nabla\F-\nabla E|^2)dxds\\
&\le -\int_0^t\int_\O(\F-E)^\top(\F-E):\nabla \u dxds
-\int_0^t\int_\O(v-\u)\cdot\nabla \u\cdot v dxds \\
&\quad -\int_0^t\int_\O(v-\u)\cdot\nabla \F: E dxds
+\int_0^t\int_\O \F:E\nabla(v-\u) dxds \\
&\quad -\int_0^t\int_\O \F^\top \F:\nabla(v-\u) dxds\\
&=-\int_0^t\int_\O(\F-E)^\top(\F-E):\nabla \u dxds
-\int_0^t\int_\O(v-\u)\cdot\nabla \u\cdot(v-\u)
dxds\\&\quad-\int_0^t\int_\O(v-\u)\cdot\nabla \F:(E-\F)
dxds+\int_0^t\int_\O (E-\F)^\top\F:\nabla(v-\u) dxds\\
&:= I,
\end{split}
\end{equation}
where, we used twice the fact
$$\int_\O v\cdot\nabla \u\cdot \u dx=0,$$
if $\Dv v=0$.
For $I$, we have, by H\"{o}lder's inequality,
\begin{equation}\label{35}
\begin{split}
|I|&\le \int_0^t (\|\nabla \u\|_{L^\infty(\O)}+\|\nabla
\F\|_{L^\infty(\O)})\left(\int_\O(|\F-E|^2+|\u-v|^2)dx\right)ds\\
&\quad+\f{1}{2}\int_0^t\int_\O|\nabla
v-\nabla\u|^2dxds+C\int_0^t\|\F\|_{L^\infty}^2\|E-\F\|^2_{L^2} ds.
\end{split}
\end{equation}
Substituting \eqref{35} back to \eqref{34}, one has
\begin{equation}\label{36}
\begin{split}
&\f{1}{2}\int_\O(|\u(t)-v(t)|^2+|\F(t)-E(t)|^2)dx+\f{1}{2}\int_0^t\int_\O(|\nabla\u-\nabla
v|^2+|\nabla\F-\nabla E|^2)dxds\\
&\le\int_0^t (\|\nabla
\u\|_{L^\infty(\O)}+\|\nabla \F\|_{L^\infty(\O)}
+C\|\F\|^2_{L^\infty(\O)})\left(\int_\O(|\F-E|^2+|\u-v|^2)dx\right)ds.
\end{split}
\end{equation}
Notice that
$$\|\nabla
\u\|_{L^\infty(\O)}+\|\nabla \F\|_{L^\infty(\O)}
+\|\F\|^2_{L^\infty(\O)}\in L^1(0,T).$$
Therefore, using \eqref{36} together with Gr\"{o}nwall's inequality, we finally conclude that $\u=v$,
$\F=E$ $a.e$ and thus $P=\Pi
$ in $\O\times(0,T)$.

The proof of Theorem \ref{T2} is complete.

\bigskip\bigskip

\section*{Acknowledgments}

Xianpeng Hu's research was supported in part by the National Science
Foundation grant DMS-0604362 and by the Mellon Predoctoral Fellowship of the University of Pittsburgh.
Dehua Wang's research was supported in part by the National Science
Foundation under Grant DMS-0604362, and by the Office of Naval Research under Grant N00014-07-1-0668.

\end{document}